\documentstyle[12pt]{article}
\author{
        {\bf Ladislav HLAVAT\'Y\thanks{e-mail  :hlavaty@br.fjfi.cvut.cz} }\\
        {\it  Department of Physics, Faculty of Nuclear Science,}\\
        {\it Czech Technical University }\\
        {\it Brehova 7, 115 19 Prague 1, Czech Republic}\\
        }

\title{Yang--Baxter systems, solutions and applications}

\def \be {\begin{equation}}
\def \bea {\begin{eqnarray}}
\def \ee {\end{equation}}
\def \eea {\end{eqnarray}}

\def \ll {\label}

\def \ox {\otimes}
\def \unit {{\bf 1}}
\def \complex {{\bf C}}

\def \newblock {}
\def \rf {(\ref}
\def \eqn {equation}
\def \sln {solution}
\def \YB {Yang--Baxter}
\def \YBE {Yang--Baxter equation}
\def \YBS {Yang--Baxter system}
\def \cfn {classification}
\def \reln {relation}
\def \tfn {transformation}
\def \mt {matri}

\begin{document}
\maketitle
q-alg/9711027
\begin{abstract}
Two types of \YBS s play roles in the theoretical physics -- constant and 
colour dependent. The constant systems are used mainly for 
construction of special Hopf algebra while the colour or spectral 
dependent for construction of quantum integrable models.
Examples of both types together with their particular solutions are 
presented. The complete \sln {} is known only for the constant system 
related to the quantized braided groups in the dimension two. 
The strategy for solution of the system related to quantum doubles is 
suggested and partial results are presented.
\end{abstract}
\section{Introduction}

The Yang--Baxter equations proved to be an important tool for various 
branches of theoretical physics. There are several types of the  
Yang--Baxter equations.

The simplest are the {\em constant}  Yang--Baxter equations. They 
represent a system of $N^6$ cubic \eqn s for elements of $N^2\times 
N^2$ matrix $R$ and can be written in the well known form
\begin{equation} R_{12}R_{13}R_{23}     =
R_{23}R_{13}R_{12}.
\ll{cybe}\end{equation}
A more complicated version of the Yang--Baxter equations is the case 
where the matrix $R$ depends on one or two parameters and the
\eqn s acquire the forms
\begin{equation} R_{12}(u)R_{13}(u+v)R_{23}(v)     =
R_{23}(v)R_{13}(u+v)R_{12}(u)
\ll{sdybe}\end{equation}
or
\begin{equation} R_{12}(u_1,u_2)R_{13}(u_1,u_3)R_{23}(u_2,u_3) =
R_{23}(u_2,u_3)R_{13}(u_1,u_3)R_{12}(u_1,u_2).
\ll{cdybe}\end{equation}
They are called {\em spectral} and {\em colour dependent} \YBE s.
It is easy to check that \rf{cybe}) and \rf{sdybe}) are special cases of 
\rf{cdybe}) for $R(u,v)=R(u-v)$ and $R(u,v)=R$ respectively.

Even though many \sln s are known for all types of the \YBE s in various 
dimensions \cite{sogo:ybecfn,hla:usybe,hla:baxtn,cinani},  the complete 
solution is known only for the constant \YBE s in the dimension two, i.e. 
matrices $4\times 4$, until now \cite{hie:ybecfn}.

Few years ago, extensions of the \YBE s for several matrices, called 
\YBS s,
appeared in literature. The goal of the present paper is to give a review 
of such systems and show results on  \cfn {} of \sln s of some constant 
systems.
\section{Review of  Yang--Baxter systems}
As the \YBS s usually contain several \YB--type \eqn s it is convenient to 
introduce the following notation: 
{\em Constant Yang--Baxter commutator} $[R,S,T]$ of (constant) 
$N^2\times N^2$ matrices $R,S,T$ is $N^3\times N^3$ matrix
\be [R,S,T]:=R_{12}S_{13}T_{23}-
T_{23}S_{13}R_{12} \ll{cybc} \ee
and
{\em spectral or colour dependent Yang--Baxter commutator} 
$[[R,S,T]]$ of (spectral or colour dependent) $N^2\times N^2$ matrices 
$R,S,T$ is $N^3\times N^3$ matrix
\[ [[R,S,T]]=[[R,S,T]] (u_1,u_2,u_3):=\]
\be R_{12}(u_1,u_2)S_{13}(u_1,u_3)
T_{23}(u_2,u_3)-
T_{23}(u_2,u_3)S_{13}(u_1,u_3)
R_{12}(u_1,u_2)  \ll{scdybc}          \ee
In this notation, the constant and spectral or colour dependent \YBE s 
read
\be [R,R,R]=0,\ [[R,R,R]]=0. \ee
\subsection{Constant \YBS s}
\subsubsection{System for quantized braided groups}
The quantized braided groups were introduced recently in
\cite{hla:qbg}  
combining Majid's concept of braided groups \cite{maj:bg} and the FRT
formulation of quantum supergroups \cite{Liao}. The
 generators $T_i^j, \  i,j \in
\{1,\ldots ,d=dimV\}$ of quantized braided groups  
satisfy the algebraic and braid relations
\begin{equation}
Q_{12}R_{12} ^{-1}T_{1}R_{12}T_{2}     =   R_{21}   ^{-1}
T_{2}R_{21}T_{1} Q_{12}             \label{qbgreln}
\end{equation}
\begin{equation}
\psi(T_{1}\ox R_{12}T_{2}  )   =   R_{12}T_{2}\ox R_{12}^{-1}
T_{1} R_{12}
\end{equation}
where the numerical matrices $Q,\ R$ satisfy the system of
Yang--Baxter--type equations
\begin{equation} [Q,Q,Q]=0,\ [R,R,R]=0,\label{rrr0}
\end{equation}
\begin{equation}
[Q,R,R]=0,\ [R,R,Q]=0. \label{zzr0}
\end{equation}

For \cfn{} of  the quantized braided groups in the dimension two
 we have to find all solutions of this Yang--Baxter system that are both 
invertible and
have the so called second inversion $(R^{t_1})^{-1}$.

The system (\ref{rrr0}) -- (\ref{zzr0}) is invariant under
\be Q'= \lambda(S\ox S)Q(S\ox S)^{-1},
 \ R'=\kappa(S\ox S)R(S\ox S)^{-1},\ \lambda,\kappa\in \complex,\
S\in SL(2,\complex) \ll{sym1}\ee
and
\be Q''=Q^+=PQP,\  R''=R^+=PRP. \ll{sym2} \ee
These symmetries are  important for \cfn{} of \sln s.

The complete set of invertible solutions in the dimension two, i.e. for 
matrices $4\times 4$, is given in \cite{hla:cfnqbg}.
Simple solutions are  $(R,\,Q=R)$  or  $(R,\,Q=PR^{-1}P),\ R$
being  solution  of  the  \YBE{} \rf{cybe}). On the other hand these \sln s 
are the only ones that solve
\rf{rrr0})--\rf{zzr0})
for all but  four classes of invertible \YB{} solutions $R$
of the Hietarinta's list \cite{hie:ybecfn}.

The exceptional cases of $R$, for which other solutions exist are
\begin{equation}
R=\left( \begin{array}{cccc}
1&0&0&0\\0&1&0&0\\0&0&1&0\\1&0&0&-1       \end{array}
\right), \
\left( \begin{array}{cccc}
0&0&0&1\\0&0&t&0\\0&t&0&0\\1&0&0&0       \end{array}
\right),\
\left( \begin{array}{cccc}
1&0&0&0\\x&1&0&0\\y&0&1&0\\z&y&x&1       \end{array}
\right), \ \ee
and diagonal matrix.
The corresponding solutions of \rf{rrr0})--\rf{zzr0})
were found in \cite{hla:cfnqbg}.

\subsubsection{System for quantum doubles}
Quantum doubles are special quasitriangular Hopf algebras 
constructed from the tensor product of Hopf algebras by defining a 
pairing between them.
In the paper  \cite{vlad:qdouble} a method of obtaining the quantum 
doubles for  pairs of FRT quantum groups is presented. 
Let two quantum groups are given by relations \cite{FRT}
\[ W_{12}U_1U_2=U_2U_1W_{12} \]
\[ Z_{12}T_1T_2=T_2T_1Z_{12} \]
where $W$ and $Z$ are matrices satisfying the \YBE s
\begin{equation} [W,W,W]=0,\ [Z,Z,Z]=0,
\ll{vlad1}\end{equation}
and suppose that there is a \mt x $X$ that satisfies the \eqn s
\begin{equation} [W,X,X]=0,\ [X,X,Z]=0.
\ll{vlad2}\end{equation}
Then the \reln s
\be X_{12}U_1T_2=T_2U_1X_{12} \ll{wxzs}\ee
define  quantum double with the pairing
\be <U_1,T_2>=X_{12}. \ee
The problem of solution of the system \rf{vlad1}),\rf{vlad2}) will be 
attacked in the
section \ref{sec:swxzs}.
\subsubsection{System for generalized reflection algebras}
The system \rf{vlad1}),\rf{vlad2}) is a special case of a more general 
\YBS{}
\[ [A,A,A]= 0, \hskip 1cm [D,D,D]=0, \]
\be [A,C,C]= 0, \hskip 1cm [D,B,B]=0, \ll{conc}\ee
\[ [A,B^+,B^+]= 0, \hskip 1cm [D,C^+,C^+]=0, \]
\[ [A,C,B^+]= 0, \hskip 1cm [D,B,C^+]=0, \]
where the superscript $X^+$ means $PXP$, and $P$ is the
{\em permutation matrix} $P_{ij}^{kl}=\delta_i^l\delta_j^k$.

This system has appeared as the consistency conditions for the 
algebra  generated by elements
$L_j^k,\
j,k \in \{1,2,\ldots ,N\}$
satisfying quadratic relations
\be A_{12}L_1B_{12}L_2=L_2C_{12}L_1D_{12} \ll{albl1} \ee
where $L=\{L_j^k\}_{j,k = 1}^{M}$ and $A,B,C,D$ are numerical matrices
$N^2\times N^2$ (i.e.
$A_{12}=\{(A)_{i_1i_2}^{j_1j_2}\in \complex \}
_{i_1,i_2,j_1,j_2 = 1}^{N}$
and similarly for $B,\ C,\ D$). These algebras were considered
in \cite{hla:afqnulm} and include (algebras of functions on) quantum
groups, quantum supergroups, braided groups, quantized braided
groups, reflection algebras and others.

\subsection{Spectral or colour dependent \YBS s}
In quantized ultralocal models, i.e. such that the Poisson brackets of 
the fields is proportional to $\delta$--function, commutation relations 
for  elements of the Lax operator  can be written as
\be R_{12}(\lambda-\mu)L_1(\lambda)L_2(\mu) = 
L_2(\mu)L_1(\lambda)R_{12}(\lambda-\mu).\ee
where the matrix $R(\lambda-\mu)$
satisfies the spectrally dependent
\YBE{} \rf{sdybe}).

The consistency conditions for the algebras used for quantization of 
nonultralocal models yield spectral or colour dependent \YBS s.

\subsubsection{Quadratic algebras for nonultralocal models}
In the paper \cite{frimai}, Freidel and Maillet suggested the 
commutation relations for the elements of the quantized Lax operator 
$L(\lambda)$ in the form
\be A_{12}(u_1,u_2)L_1(u_1)B_{12}(u_1,u_2)L_2(u_2)= 
L_2(u_2)C_{12}(u_1,u_2)L_1(u_1)D_{12}(u_1,u_2) \ll{albl} \ee
where  $A,B,C,D$ are spectral or colour dependent $N^2\times N^2$ 
matrices 
that satisfy the \eqn s of the \YBS{}
\[ [[A,A,A]]= 0, \hskip 1cm [[D,D,D]]=0, \]
\be [[A,C,C]]= 0, \hskip 1cm [[D,B,B]]=0, \ll{scdconc}\ee
\[ [[A,B^\ddagger,B^\ddagger]]= 0, \hskip 1cm 
[[D,C^\ddagger,C^\ddagger]]=0, \]
\[ [[A,C,B^\ddagger]]= 0, \hskip 1cm [[D,B,C^\ddagger]]=0. \]
The superscript $X^\ddagger$ is defined by $X^\ddagger(u,v):=PX(v,u)P$.
This algebra with $C=B^\ddagger$ was introduced also in 
\cite{nijcap:lnp,nijetal:pra92}.

Several particular  solutions are known,
e.g.
\[ A(u,v)=(u-v)\unit+P \]
\be B^\ddagger(u,v)=C(u,v)=v\unit+\sigma_-\ox\sigma_+\ee
\[ D(u,v)=u-v+(1-u/v)\sigma_-\ox\sigma_+(1-v/u)\sigma_+\ox\sigma_- \]
where 
\[sigma_+=\left( \begin{array}{cc}
0&1\\0&0 \end{array}\right),\ \ \sigma_-=\left( \begin{array}{cc}
0&0\\1&0 \end{array}\right). \]

\subsubsection{Multiple braided product of quadratic 
algebras}\label{mpa}
The crucial object
for construction of quantum integrable models is
the monodromy matrix \cite{takhfad:UNM}. In the quantized models it is
a representation of the multiple
matrix coproduct $\Delta^N(L)$ for a quadratic algebra generated by 
$L$. The usual matrix
coproduct cannot be used for nonultralocal models unless we 
introduce the braiding
structure to the tensor products of the algebras \cite{hla:afqnulm}. The 
structure can be expressed as the algebra
generated by
$N\times M^2$ generators
\be (L^{I})_j^k(u),\ I \in \{1,2,\ldots ,N\},\
j,k \in \{1,2,\ldots ,M\} \ll{lgen}\ee
satisfying quadratic relations
\be
L_1^J(u_1)X_{12}^{JK}(u_1,u_2)L_2^K(u_2)=
W_{12}^{JK}(u_1,u_2)L_2^K(u_2)
Y_{12}^{JK}(u_1,u_2)L_1^J(u_1)
Z_{12}^{JK}(u_1,u_2)
\ll{ljxlklam} \ee
where $X^{JK},\ W^{JK},\ Y^{JK},\ Z^{JK}$ for fixed
$J,K\in \{1,2,\ldots ,N\}$ are numerical 
invertible colour dependent  $M^2\times M^2$ matrices. 
No summation over indices
$I,J,K,\ldots$ is assumed.

The importance of these algebra consists in the fact that we can
find a  commuting subalgebra that can be used for construction of
quantum hamiltonian of a model together with conserved
quantities.

The  quadratic algebras defined by relations of the form
(\ref{ljxlklam}) must satisfy
consistency conditions of the Yang--Baxter--type
 that follow from the requirement that no
supplementary higher degree relations are necessary for unique
transpositions of three and more elements \cite{FRT}. For relations of
the form
(\ref{ljxlklam})  the conditions read (cf.\rf{scdconc}))
\be \{ [Z,Z,Z]\} =0,\ \{ [W,W,W]\} =0 \ll {sybe1}\ee
\be \{ [Z,X,X]\} =0,\ \{ [X,X,W]\} =0 \ll {sybe2}\ee
\be \{ [Z,Y^{\ddagger},Y^{\ddagger}]\} =0,\
\{ [Y^{\ddagger},Y^{\ddagger},W]\} =0
\ll {sybe3}\ee
\be \{ [Z,X,Y^{\ddagger}]\} =0,\ \{ [Y^{\ddagger},X,W]\} =0 \ll {sybe4} \ee
where  by $\{ [R,S,T]\} =0$ we mean that
\[  [[R^{J_1J_2}S^{J_1J_3}T^{J_2J_3}]] =
R_{12}^{J_1J_2}(u_1,u_2)S^{J_1J_3}_{13}(u_1,u_3)T^{J_2J_3}_{23}(u_2,u_
3)\] \be
-T_{23}^{J_2J_3}(u_2,u_3)S^{J_1J_3}_{13}(u_1,u_3)R^{J_1J_2}_{12}(u_1,u_
2)=0
 \ll{sybc} \ee
for all $ J_1,J_2,J_3 \in \{1,\ldots,N\}\ \}$ and 
\[(Y^{\ddagger})^{JK}(u,v):=(Y^{KJ}(v,u))^+ = PY^{KJ}(v,u)P. \]

Several particular solutions of \rf{sybe1}) -- \rf{sybe4}) are
given in \cite{hla:afqnulm}.

\section{Solving the \YBS {} for the quantum double}\ll{sec:swxzs}
Our longtime goal is solution of the spectral dependent \YBS s for 
nonultralocal models presented in the previous section. Nevertheless 
to get a deeper experience with such complicated systems we decided 
to solve their constant  versions first. One can easily check that the 
constant version of the system \rf{scdconc}) where $B=C^\ddagger$ 
can be written as the system for the quantum doubles 
\begin{equation} [W,W,W]=0,\ll{www}\ee
\begin{equation} [W,X,X]=0,\ll{wxx}\ee
\be [X,X,Z]=0, \ll{xxz}\end{equation}
\be [Z,Z,Z]=0, \ll{zzz}\ee

We shall shall look for  \sln s of this system in the lowest nontrivial 
dimension two (i.e. matrices $4\times 4$).

Solution of the system is essentially facilitated by knowledge of 
symmetries of the system \rf{www})--\rf{zzz}).
The space of \sln s is invariant under both continuous transformations
\be W'=\omega(T\ox T)W(T\ox T)^{-1}\ll{wsym}\ee
\be X'=\xi(T\ox S)X(T\ox S)^{-1}\ll{xsym}\ee
\be Z'=\zeta(S\ox S)Z(S\ox S)^{-1}\ll{zsym}\ee
where
\[ \omega,\xi,\zeta\in\complex,\ \ T,S\in SL(2,\complex) \]
and the 
discrete transformations
\be (W',X',Z')=(W^t,X^t,Z^t) \ll{dsym}\ee
\be (W',X',Z')=(W^a,X,Z^b), \ \ a=id,\# ,\ b=id,\#  \ll{dsym1}\ee
\be (W',X',Z')=(W^c,X^-,Z^d), \ \ c=+,-,\ d=+,- \ll{dsym2}\ee
\be (W',X',Z')=(Z^c,X^+,W^d), \ \ c=+,-,\ d=+,- \ll{dsym3}\ee
where $Y^t$ is transpose of $Y$, $Y^+:=PYP,\ Y^-:=Y^{-1},$ $Y^{id}:=Y,\ Y^\# 
:=(Y^+)^-=(Y^-)^+$.

Beside that one can guess several simple solutions of the system, 
\rf{www})--\rf{zzz}),namely
\be W=X=Z=R,\ {\rm where} \ [R,R,R]=0, \ee
\be X,=\unit,\ W,Z\ {\rm arbitrary\ solutions\ of}\  [W,W,W]=0,\ [Z,Z,Z]=0, \ee
and
\be W=Z=P,\ X\ {\rm arbitrary\ \mt x},\ee
where $P$ is the permutation matrix and $X$ is arbitrary.

The last two solutions give positive answers to the following questions.
Is there a \mt x $X$ such that for any pair of matrices $W,Z$ that solve 
the \YBE s  the triple $(W,X,Z)$ solves the system \rf{www})--\rf{zzz})?
Is there for any \mt x $X$ a pair of matrices $W,Z$  such that the triple 
$(W,X,Z)$ solves the system \rf{www})--\rf{zzz})?

Besides the above mentioned, one can find  solutions of the system 
\rf{www})--\rf{zzz}) from the knowledge of solutions of the system 
(\ref{rrr0}--\ref{zzr0}). Namely, if $(Q,\, R)$ is a solution of
the system of Yang--Baxter type equations (\ref{rrr0}--\ref{zzr0})
then $(W=Q,X=R,Z=R^+QR^-)$ is a solution of the
system \rf{www})--\rf{zzz}). 

The question is if there are other \sln s than those given above. The 
answer is also positive and the author believes that the complete \sln {} 
of the system \rf{www})--\rf{zzz}) in the dimension two can be found.
\subsection{Strategies of solution}
Even for the lowest nontrivial dimension two, the system 
\rf{www})--\rf{zzz}) represents a tremendous task - solving 256 cubic 
equations for 48 unknowns. That's why it is understandable that the 
assistance of computer programs for symbolic calculations is essential 
in the following\footnote{Reduce 3.6 was used}. On the other hand it 
does not mean that one can find the solutions by pure brute force, 
namely applying a procedure SOLVE to the system of the 256 
equations.

There are several strategies for solution of the above given problem. 
All of them are based on the knowledge of the symmetries of the 
system and the knowledge of the complete set of solutions of the 
\YBE{} in the dimension two.

One possible (and obvious) strategy is solving the \eqn s 
\rf{wxx}),\rf{xxz}) for all pairs of matrices $W,Z$ in the Hietarinta's list of 
\sln s of the \YBE{} \cite{hie:ybecfn}. By this way we reduce the problem 
to 128 quadratic \eqn s for 16-22 unknowns (depending on the number 
of parameters in the \sln s of the \YBE{}).

Another strategy is to use the symmetry \rf{xsym}) to simplify the matrix 
$X$ as much as possible, then solve the linear \eqn s \rf{wxx}),\rf{xxz}) 
for $W$ and $Z$ and finally solve the \YBE s \rf{www}),\rf{zzz}). This is 
an analogue of the strategy accepted in \cite{hla:cfnqbg} for solving the 
system \rf{rrr0})--\rf{zzr0}).

Both the mentioned strategies yield sets of equations and unknowns 
that are still too large to be solved by  the computer program.
The strategy that seems to be working is a combination of the two 
previous. We start with a solution of the \YBE{} $W$ and solve $X$ from 
the \eqn{} \rf{wxx}) using the symmetries \rf{wsym})--\rf{dsym3}) that 
leave $W$ form--invariant (i.e. changing only values of its parameters). 
Then we solve the \eqn{}  \rf{xxz}) (linear in elements of $Z$)  and  
determine $Z$ from \rf{zzz}) using the results of the previous step. The 
results of this strategy for the so called (non)standard \sln{} of the 
\YBE{} are given in the next subsection.
\subsection{Results for the standard and nonstandard \sln{} of the 
\YBE}
Let 
\begin{equation}
W=W_{q,s}=\left( \begin{array}{cccc}
q&0&0&0\\0&s^{-1}&0&0\\0&q-q^{-1}&s&0\\0&0&0&t       \end{array}
\right),\ t=q\ {\rm or}\ t=-q^{-1},\ q,s\in\complex\setminus\{0\},\ q^2\neq 1. 
\ll{w5}\ee
This \mt x is called standard, for $t=q$, or nonstandard, for $t=-q^{-1}$, 
\sln{} of the \YBE{}.
Up to the transformations \rf{wsym})--\rf{dsym3}) we get the following 
list of solutions
of the \eqn{} \rf{wxx}) for any (invertible) $q,s\in\complex$. 
\begin{equation}
X_1=\left( \begin{array}{cccc}
a&0&0&0\\c&a&0&0\\0&0&b&0\\0&0&d&b       \end{array}
\right), 
\ \ X_2=\left( \begin{array}{cccc}
q&0&0&0\\0&s^{-1}&0&0\\0&a&b&0\\0&0&0&\frac{b}{s}t       \end{array}
\right), \ \ 
\ll{w5x2}\ee
\begin{equation}
X_3=\left( \begin{array}{cccc}
a&0&0&0\\0&b&0&0\\0&0&c&0\\0&0&0&d      \end{array}
\right). \ll{w5x4}\ee
Beside that  we get special solutions 
\be X_4=\left( \begin{array}{cccc}
q&0&0&c\\0&s^{-1}&0&0\\0&a&b&0\\0&0&0&\frac{b}{s}t       \end{array}
\right), \ \ {\rm for}\ s^2=1. \ee
\begin{equation}
X_5=\left( \begin{array}{cccc}
a&0&0&b\\0&-a&b&0\\0&0&c&0\\0&0&0&c       \end{array}
\right),\ \ \ {\rm for}\ q=-s=i.\ee
\be X_6=\left( \begin{array}{cccc}
0&0&ia&0\\2iab/c&0&0&a\\ib&0&0&c\\0&b&0&0       \end{array}
\right), \ \ {\rm for}\ q=i,\ s=1.\ll{w5x5}\ee
respectively.

Each of these \sln s of the \eqn {} \rf{wxx}) can be extended to the \sln{} 
of the system \rf{www})--\rf{zzz}) by $Z=P$ where $P$ is the 
permutation matrix
\begin{equation}
P=\left( \begin{array}{cccc}
1&0&0&0\\0&0&1&0\\0&1&0&0\\0&0&0&1     \end{array}
\right).\ee 
For $X_1$
--$X_5$ there are still other solutions.

Solutions of the \eqn s  \rf{xxz})--\rf{zzz}) where $X=X_1$ are
\begin{equation}
Z_{10}=\left( \begin{array}{cccc}
1&0&0&0\\x&1&0&0\\y&0&1&0\\z&y&x&1       \end{array}
\right),\ \ 
Z_{11}=\left( \begin{array}{cccc}
1&0&0&0\\x&1&0&0\\-x&0&1&0\\-xy&-y&y&1       \end{array}
\right),\ee
so that the  triples
$(W_{q,s},X_1,Z_{11})$ and $(W_{q,s},X_1,Z_{12})$ are solutions to the 
system \rf{www})--\rf{zzz}).

Similarly, for $X_2$ we get the \sln s of the \eqn s  \rf{xxz})--\rf{zzz})  
\begin{equation}
Z_{20}=\left( \begin{array}{cccc}
q&0&0&0\\0&b^{-1}&0&0\\0&q-q^{-1}&b&0\\0&0&0&t \end{array}
\right). \ll{w5x3}\ee
and 
if $q^2=-1$ then moreover
\begin{equation}
Z_{21}=\left( \begin{array}{cccc}
q&0&0&\delta\\0&r&0&0\\0&q-rbq^{-1}&b&0\\0&0&0&-rbq^{-1}\end{array}
\right), \ee
where $\delta=0$ if $b^2\neq -1$.

For the diagonal $X_3$ we find that any "six-vertex" matrix Z solves 
\rf{xxz}) so that the solutions of \rf{xxz})--\rf{zzz}) are 
\be Z_{30}=\left( \begin{array}{cccc}
p&0&0&0\\0&r&0&0\\0&0&x&0\\0&0&0&y       \end{array}
\right), \ee
\begin{equation}
Z_{31}=\left( \begin{array}{cccc}
p&0&0&0\\0&r^{-1}&0&0\\0&p-p^{-1}&r&0\\0&0&0&p       \end{array}
\right), \ 
Z_{32}=\left( \begin{array}{cccc}
p&0&0&0\\0&r^{-1}&0&0\\0&p-p^{-1}&r&0\\0&0&0&-p^{-1}       \end{array}
\right).\ll{z4}\ee

If $a=b,c=-d$ in $X_3$ then any "eight-vertex" matrix Z solves \rf{xxz}). 
To obtain all invertible \sln s of the system \rf{xxz})--\rf{zzz}) in this case 
one must perform the \cfn{} of the "eight-vertex" \sln s of the \YBE 
\rf{zzz}) up to the symmetries \rf{zsym}) where $S$ is the (anti)diagonal 
\mt x because only \tfn s \rf{xsym}) with (anti)diagonal $S$   leave $X_3$ 
form--invariant. This \cfn {} was done in \cite{hla:cfnqbg} and beside the 
"eight-vertex" \sln s given in \cite{hla:usybe} contains the \mt ces
\be Z=\left( \begin{array}{cccc}
x&0&0&y\\0&\pm x&y&0\\0&y&\pm x&0\\y&0&0&x      \end{array}
\right),\ \ x,y\neq 0                             \ll{spec8v}\ee

If $a=b,c=d$ then any triple
$(W_{q,s},X_3,Z)$ where $Z$ solves the \YBE{}, is solution to the system 
\rf{www})--\rf{zzz}). The list of these solutions is then equivalent to the 
list of \sln s in \cite{hie:ybecfn}.

The only \mt x $Z$ that solves the \eqn s \rf{xxz}) for $X=X_4$ where $q^2 
\neq -1$ or $b^2 \neq 1$ is $Z=P$. If $q^2=-b^2=-1$ then there is another 
solution, namely
\be 
Z_{41}=\left( \begin{array}{cccc}
p&0&0&\frac{ac}{2}(p+p^{-1})\\0&bp^{-1}&0&0\\0&p-p^{-1}&bp&0\\
0&0&0&-p^{-1}       \end{array}
\right).\ll{z41}\ee

Finally for  $q=-s=i$ we get solutions $(W_{i,-i},X_5,Z)$ where
\be Z=Z_{51}=\left( \begin{array}{cccc}
1&0&0&1\\0&\epsilon&1&0\\0&1&-\epsilon&0\\-1&0&0&1     \end{array}
\right),\ee
or
\be Z=Z_{52}=\left( \begin{array}{cccc}
k-k^{-1}+2&0&0&k-k^{-1}\\0&k+k^{-1}&k-k^{-1}&0\\ 
0&k-k^{-1}&k+k^{-1}&0\\k-k^{-1}&0&0&k-k^{-1}-2     \end{array}
\right), 
\ee 
or 
\be 
Z=Z_{53}=\left( \begin{array}{cccc}
k&0&0&0\\0&\epsilon k&0&0\\0&k-1&1&0\\
\epsilon (k-1)&0&0&-1       \end{array}
\right),\ee
or
\be 
Z=Z_{54}=\left( \begin{array}{cccc}
k&0&0&0\\0&1&0&0\\0&k-k^{-1}&1&0\\
0&0&0&-k^{-1}      \end{array}
\right),\ee
where $ \epsilon^2=1,\ k=c/a$.
\section{Conclusions}
Recently two types of \YBS s  have appeared in the theoretical physics 
-- constant and colour dependent. The constant systems are used 
mainly for construction of special Hopf algebras while the colour or 
spectral dependent for construction of quantum integrable models.
We have presented examples of both types together with their 
particular solutions. The complete \sln {} is known only for the constant 
system \rf{rrr0})--\rf{zzr0}) in the dimension two. 

The strategy for solution of the quantum double system 
\rf{vlad1}),\rf{vlad2}) was suggested and partial results were presented.
Namely, if $W$ is the standard or nonstandard \sln{} \rf{w5}) of the 
\YBE{} then the complete set of invertible \sln s of this system up to 
symmetry \tfn s \rf{wsym})--\rf{dsym3}) is formed by the triples
\[ (W_{q,s},X,P) \ {\rm where} \ X=X_1,X_2,X_3 \]
\[ (W_{q,s},X_1,Z_{11}),\ (W_{q,s},X_1,Z_{12}) \] 
\[ (W_{q,s},X_2,Z_{20}),\ (W_{i,s},X_2(q=i),Z_{21})  \] 
\[ (W_{q,s},X_3,Z_{30}),\ (W_{q,s},X_3,Z_{31}),\ 
(W_{q,s},X_3,Z_{32}),\ 
  (W_{q,s},diag(a,-a,b,b),Z_{8V}),\]
where $Z_{8V}$ is invertible "eight-vertex" \sln{} of the \YBE {} (for the 
\cfn{} see \cite{hla:usybe, hla:cfnqbg}), 
\[  (W_{q,s},diag(a,a,b,b),Z),\]
where $Z$ is arbitrary invertible \sln{} of the \YBE {} (for the \cfn{} see 
\cite{hie:ybecfn}), 
\[ (W_{q,\pm 1},X_4,P),\ (W_{i,\pm 1},X_4(q=i),Z_{41}), \]
\[ (W_{i,-i},X_5,P),\ (W_{i,-i},X_5,Z_{51}),\ (W_{i,-i},X_5,Z_{52}),\ 
(W_{i,-i},X_5,Z_{53}),\ (W_{i,-i},X_5,Z_{54}),\ \]
and $ (W_{i,1},X_6,P)$.
The forms of the \mt ces
$P$, $X_i$ and $Z_{ij}$ are given in the previous section.

From this result it is clear that quite distinct types of FRT quantum 
groups can be combined into the quantum doubles.

\end{document}